\documentclass{amsart}
\usepackage{amssymb, amsfonts}
\usepackage[v2]{xy}

\newtheorem{thm}{Theorem}[section]

\newtheorem{prop}[thm]{Proposition}
\newtheorem{lem}[thm]{Lemma}

\theoremstyle{definition}
\newtheorem{defn}[thm]{Definition}

\newtheorem{con}[thm]{Construction}

\theoremstyle{remark}
\newtheorem{rem}[thm]{Remark}

\DeclareFontFamily{OMS}{rsfs}{\skewchar\font'60}
\DeclareFontShape{OMS}{rsfs}{m}{n}{<-5>rsfs5 <5-7>rsfs7 <7->rsfs10 }{}
\DeclareSymbolFont{rsfs}{OMS}{rsfs}{m}{n}
\DeclareSymbolFontAlphabet{\scr}{rsfs}

\let\overto\xrightarrow

\newcommand{\sB}{\scr{B}}
\newcommand{\sC}{\scr{C}}
\newcommand{\sD}{\scr{D}}

\newcommand{\sU}{\scr{U}}
\newcommand{\sV}{\scr{V}}

\newcommand{\epz}{\varepsilon}
\newcommand{\ph}{\phi}

\newcommand{\la}{\lambda}
\newcommand{\tha}{\theta}

\newcommand{\si}{\sigma}
\newcommand{\ta}{\tau}

\newcommand{\ze}{\zeta}
\newcommand{\om}{\omega}

\newcommand{\SI}{\Sigma}

\newcommand{\com}{\circ}     % composition of functions
\newcommand{\iso}{\cong}     % preferred isomorphism symbol
\newcommand{\htp}{\simeq}    % homotopy symbol
\newcommand{\sma}{\wedge}    % smash product 
\newcommand{\rtarr}{\longrightarrow}

\def\quickop#1{\expandafter\newcommand\csname #1\endcsname{\operatorname{#1}}}
\quickop{Hom} \quickop{End} \quickop{Aut} \quickop{Tel} \quickop{Mic} 
\quickop{Ext} \quickop{Tor} \quickop{Id} \quickop{Coker} \quickop{Ker}
\quickop{Lim} \quickop{Colim} \quickop{Holim} \quickop{Hocolim}
\quickop{id} \quickop{tel} \quickop{mic} \quickop{coker} \quickop{Map}
\quickop{colim} \quickop{holim} \quickop{hocolim} \quickop{im}

\makeatletter
\let\c@equation\c@thm
\makeatother
\numberwithin{equation}{section}

\let\overto\xrightarrow

\title{The Wirthm\"uller isomorphism revisited}

\author{J.P. May}

\thanks{The author was partially supported by the NSF}	

\begin{document}

\begin{abstract}
We show how the formal Wirthm\"uller isomorphism theorem proven in \cite{May}
simplifies the proof of the Wirthm\"uller isomorphism in equivariant stable 
homotopy theory. Other examples from equivariant stable homotopy theory show 
that the hypotheses of the formal Wirthm\"uller and formal Grothendieck
isomorphism theorems in \cite{May} cannot be weakened.
\end{abstract}

\maketitle

\tableofcontents

We illustrate the force of the formal Wirthm\"uller isomorphism theorem
of \cite{May} by giving a worked example of an interesting theorem whose 
proof it simplifies, namely the Wirthm\"uller isomorphism theorem in 
equivariant stable homotopy theory. It relates categories of $G$-spectra and 
$H$-spectra for $H\subset G$. We also say just a little about the analogous 
Adams isomorphism that relates categories of $G$-spectra and $J$-spectra for a 
quotient group $J$ of $G$. That context gives an interesting situation where the 
formal hypotheses of the formal Wirthm\"uller isomorphism theorem hold but the 
conclusion fails, showing that the more substantive hypothesis is essential. 

In their general form, the Wirthm\"uller and Adams isomorphisms are due to 
Gaunce Lewis and myself \cite{LMS}. It is a pleasure to thank Lewis 
for recent e-mails, ongoing discussions, and his longstanding quest 
for simplifications and generalizations of these theorems. The analogy 
between these isomorphisms in topology and Verdier duality was first 
explored by Po Hu \cite{PoHu}, who carried out an idea of Lewis that these 
isomorphisms could be obtained using parametrized equivariant spectra. 
She obtained a substantial generalization of the Wirthm\"uller isomorphism, 
but at the price of greatly increased complexity. Exactly as in our proof 
here, the theory of \cite{May} allows considerable simplification of her work,
and that is part of our motivation. The theory of \cite{May} should also 
simplify the proof of the Wirthm\"uller isomorphism in $\mathbf{A}^1$ stable
homotopy theory that Hu proved in \cite{Hu2}. Still another example recently 
studied by Hu deals with change of universe. We show that it gives an 
interesting naturally occurring situation in which all but one of the hypotheses 
of the formal Grothendieck isomorphism of \cite{May} hold, but the conclusion 
fails because the relevant left adjoint fails to preserve compact objects.

\section{The Wirthm\"uller isomorphism}

Let $H$ be a (closed) subgroup of a compact Lie group $G$
and let $f:H\rtarr G$ be the inclusion. Let $L$ be the tangent 
$H$-representation at the coset $eH\in G/H$. Thus, if $G$ 
is finite or, more generally, if $H$ has finite index in $G$, 
then $L=0$. 

Let $\sD$ and $\sC$ be the stable homotopy categories of 
$G$-spectra and of $H$-spectra, as constructed in \cite{LMS}
or, in more modern form, \cite{MM}.  The category $\sD$
depends on a choice of a ``$G$-universe'' $\sV$ on which to 
index $G$-spectra.  We may think of $\sV$ as a
chosen collection $\{V\}$ of $G$-representations $V$ that 
contains the trivial representation and is closed under direct 
sums. The most interesting example is the complete $G$-universe
obtained by allowing all representations of $G$.  We think of 
representations as finite dimensional $G$-inner product spaces 
and let $S^V$ denote the one-point compactification of $V$. 
The point of the choice of $G$-universe is that, in the 
construction of $\sD$, we force suspension by some
representations $V$ to be equivalences $\sD\rtarr \sD$, 
and we must choose which representations to invert in this 
sense.

We insist that $G/H$ embed in a representation $V$ in our 
$G$-universe $\sV$, which is otherwise unrestricted. We index 
$H$-spectra on the $H$-universe $f^*\sV = \{f^*V\}$, where 
$f^*V$ denotes the $G$-representation $V$ viewed as an 
$H$-representation by pullback along $f$. Similarly, a $G$-spectrum $Y$ gives 
an $H$-spectrum $f^*Y$ by pullback along $f$. The functor $f^*$ has a left adjoint
$f_{!}$ and a right adjoint $f_*$. The former is usually
written as either $G_+\sma_H X$ or $G\ltimes_H X$, and the
latter is usually written as either $F_H(G_+,X)$ or
$F_H[G,X)$. The Wirthm\"uller isomorphism reads as follows.

\begin{thm}[Wirthm\"uller isomorphism]\label{Wirth} There 
is a natural isomorphism
$$\om: f_* X\rtarr f_{\sharp}X, \ \ \text{where}  \ \
f_{\sharp}X = f_{!}(X\sma S^{-L}).$$
That is, for an $H$-spectrum $X$, $F_H(G_+,X)\iso G_+\sma_H(\SI^{-L}X)$. 
\end{thm}

Here the suspension $H$-spectrum of $S^L$ is invertible 
with inverse $S^{-L}$, allowing the definition $\SI^{-L}X = X\sma S^{-L}$.
Indeed, an embedding of $G/H$ in $V\in \sV$ induces an inclusion 
$L\subset V$ of $H$-representions with 
orthogonal complement $W$. Since $S^L\sma S^W\iso S^{f^*V}$
is invertible in $\sC$, $S^L$ is also invertible in $\sC$. 

The unit objects in $\sC$ and
$\sD$ are the sphere spectra $S_H$ and $S_G$. Both
$\sC$ and $\sD$ are closed symmetric monoidal 
categories under their smash product and
function spectrum functors $\sma$ and $F$. It is 
immediate from the definitions that $f^*$ is strong
symmetric monoidal and commutes with $F$, as required
in the Wirthm\"uller context discussed in \cite[\S\S2,\,4]{May}.

\begin{rem} In this context, the projection formula would assert that
$$Y\sma F_H(G_+,X) \iso F_H(G_+,f^*Y\sma X),$$
which is false on the spectrum level but which holds on the stable category
level as a consequence of the Wirthm\"uller isomorphism. Note that, on the
spectrum level, $F_H(G_+,X)^G \iso X^H$, by a comparison of adjunctions.
\end{rem}

The isomorphism \cite[4.3]{May} required in the Wirthm\"uller
context can be written
\begin{equation}\label{At}
 D(G_+\sma_H S_H)\iso G_+\sma_H S^{-L}.
\end{equation}
It is a special case of equivariant Atiyah duality for smooth 
$G$-manifolds, which is proven by standard space level techniques 
(e.g. \cite[III\S 5]{LMS}) and is independent of the Wirthm\"uller 
isomorphism. By the tubular neighborhood theorem, we can extend an 
embedding $i: G/H\rtarr V$ to an embedding 
\begin{equation}\label{tube}
\tilde i: G\times_H W\rtarr V
\end{equation} 
of the normal bundle $G\times_H W$. Atiyah duality asserts that the 
$G$-space $G/H_+$ is Spanier-Whitehead $V$-dual to the Thom complex 
$G_+\sma_H S^W$ of the 
normal bundle of the embedding. Desuspending by $S^V$ 
in $\sD$ gives the required isomorphism (\ref{At}). 

The category $\sD$ is triangulated, with distinguished 
triangles isomorphic to canonical cofiber sequences of
$G$-spectra. The triangulation is compatible with the 
closed symmetric monoidal structure in the sense discussed 
in \cite{May2}. Moreover, $\sD$ is compactly generated. 
Writing $S^n$ for the $n$-sphere $G$-spectrum, we
can choose the generators to be the $G$-spectra $G/J_+\sma S^n$,
where $J$ ranges over the (closed) subgroups of $G$ and $n$ runs
over the integers.  The same statements apply to $\sC$,
and the functor $f^*$ is exact since it commutes with cofiber 
sequences on the level of spectra. Moreover, $f^*$ takes
generators to compact objects. Indeed, this depends only on
the fact that the $H$-spaces $G/J$ are compact and of the 
homotopy types of $H$-CW complexes, although it is easier to 
verify using the stronger fact that the $G/J$ can be decomposed
as finite $H$-CW complexes.  

\begin{rem} Observe that the generators $H/K_+$ of $\sC$ and the generators 
$G/J_+$ of $\sD$, other than $G/H_+$, need not be dualizable if the universe 
is incomplete. In fact, Lewis \cite[7.1]{Lewis} has proven that $G/J_+$ is 
dualizable if and only if $G/J$ embeds in a representation in $\sV$. 
Conceptually, it is compactness rather than dualizability of 
the generators that is relevant.
\end{rem}

Digressively, there is an interesting conceptual point to
be made about the choice of generators. There are (at least) 
three different, but Quillen equivalent, model categories of 
$G$-spectra. We can take ``$G$-spectra'' to mean $G$-spectra as 
originally defined \cite{LMS}, $S_G$-modules as defined in 
\cite{EKMM}, or orthogonal $G$-spectra as defined in \cite{MM}, 
where the cited Quillen equivalences are proven. These 
three categories are also compactly generated in the model 
theoretic sense. In the first two cases, the generators in the 
model theoretic sense can be taken to be the generators in the 
triangulated category sense that we have just specified. 
As explained in \cite{MM}, we can alternatively take all 
$G/J_+\sma S^{V-V'}$ for $V, V'\in\sV$ as generators in the
model theoretic sense. However, in the model category of orthogonal 
$G$-spectra, we not only can but must take this larger 
collection as generators in the model theoretic sense.
Nevertheless, the smaller collection suffices to generate the
associated triangulated homotopy category, since that is 
triangulated equivalent to the homotopy category obtained 
from the other two model categories.

By the formal Wirthm\"uller isomorphism theorem \cite[6.1]{May}, 
to prove Theorem \ref{Wirth}, it remains to prove that the map 
$\om:f_*X\rtarr f_{\sharp}X$ specified in \cite[4.7]{May} is an 
isomorphism when $X$ is a generating 
object $H/K_+\sma S^n$ of $\sC$. Since it is obvious that $\om$ 
commutes with suspension and desuspension, we need only consider 
the case $n=0$, where the generators in question are the suspension
$H$-spectra $\SI^{\infty}H/K_+$. 

Here is the punchline. Suppose that $G$ is finite or, more generally,
that $H$ has finite index in $G$. Then, as an $H$-space, $H/K_+$ is
a retract of the $G$-space $G/K_+$. The retraction sends cosets of
$G/K_+$ not in the image of $H/K_+$ to the disjoint basepoint of $H/K_+$.
By \cite[4.13]{May}, it follows formally that $\om$ is an isomorphism
when $X = \SI^{\infty}H/K_+$. This simple argument already completes the 
proof of Theorem \ref{Wirth} in this case.

To prove Theorem \ref{Wirth} in general, we apply \cite[4.14]{May}, 
which allows us to work one generating object at a time. This reduces 
all of our work to consideration of suspension spectra and thus of spaces.
The argument is essentially the same as part of the argument in 
\cite[III\S\S5,\,6]{LMS}, but we shall run through the space-level 
details in \S4 in order to have a readable and self-contained account.

\section{The Adams isomorphism}

Let $N$ be a normal subgroup of a compact Lie group $G$ and let $f:G\rtarr J$, $J=G/N$, 
be the quotient homomorphism. Fix a $G$-universe $\sV=\{V\}$ and index $J$-spectra on the
$N$-fixed $J$-universe $\sV^N = \{V^N\}$. Regarding $J$-representations as $G$-representations
via $f$, we obtain a second $G$-universe $f^*\sV^N =\{f^*V^N\}$, and we insist that 
$f^*\sV^N$ be contained in the original $G$-universe $\sV$. Let $\sD$ be the stable 
homotopy category of $J$-spectra indexed on $\sV^N$ and let $\sC$ be the stable 
homotopy category of $G$-spectra indexed on $f^*\sV^N$. 

Regarding $J$-spectra as $G$-spectra via $f$, we obtain a functor 
$f^*:\sD\rtarr \sC$. Its left adjoint $f_{!}$ is just the orbit spectrum functor that 
sends a $G$-spectrum $X$ indexed on $f^*\sV^N$ to $X/N$. Its right adjoint $f_*$ is just the 
fixed point spectrum functor that sends $X$ to $X^N$. The functor $f^*$ is strong symmetric
monoidal, the isomorphism \cite[2.6]{May} holds in the form $(f^*Y\sma X)/N\iso Y\sma (X/N)$,
and $f^*$ takes generating objects to compact objects. Since $G$-spectra in $\sC$ are indexed 
on an $N$-trivial universe, $S_G/N \iso S_J$ and \cite[4.3]{May} holds in the trivial form 
$D(S_J)\iso S_G/N$. Thus all of the formal hypotheses of the formal Wirthm\"uller isomorphism 
theorem, \cite[6.1]{May},  are satisfied. However, the conclusion fails, because $\om$ 
is an isomorphism on some but not all generators. Let $\sC/N$ be the thick subcategory of 
$\sC$ generated by all $\SI^n\SI^{\infty}_G G/H_+$ such that $N\subset H$.

\begin{prop} For $X\in \sC/N$, the natural map
$$\om: X^N\rtarr X/N$$
is an isomorphism of $J$-spectra.
\end{prop}
\begin{proof} 
The map $\ta:T\rtarr f_{!}C$ of \cite[\S4]{May} is here just the isomorphism $S_J\iso S_G/N$, 
and the map $\om$ of \cite[4.7]{May} is just the composite of the isomorphism $X^N\iso f^*X^N/N$ 
and $\epz/N: f^*X^N/N\rtarr X/N$. When $X=\SI^n\SI^{\infty}_G G/H_+$ with $N\subset H$, the 
latter map is also an isomorphism, but it is not an isomorphism in general.
\end{proof}

The Adams isomorphism is more subtle. Let $\sB$ be the stable homotopy category of $G$-spectra
indexed on $\sV$. We have the adjoint pair $(i_*,i^*)$ of change of universe functors 
$i_*:\sC\rtarr \sB$ and $i^*:\sB\rtarr \sC$ induced by the inclusion $i:f^*\sV^N\rtarr \sV$.
We are most interested in the composite adjunction $(i_*f^*,f_*i^*)$. Thinking of the case when 
$\sV$ is a complete $G$-universe, it is usual to regard the composite $f_*i^*:\sB\rtarr \sD$ as the 
$N$-fixed point spectrum functor from $G$-spectra indexed on $\sV$ to $J$-spectra indexed 
on $\sV^N$.

The conjugation action of $G$ on $N$ gives rise to an action of $G$ on the tangent space 
$A=A(N;G)$ of $N$ at $e$. We call $A$ the adjoint representation of $G$ at $N$. Of course,
$A=0$ if $N$ is finite. The Adams isomorphism reads as follows. Let $\sC_{f}$ be the full 
subcategory of $N$-free $G$-spectra in $\sC$. 

\begin{thm}[Adams isomorphism]\label{Adams} For $X\in \sC_{f}$, there is a natural 
isomorphism 
$$\om: f_* i^*i_*X \rtarr f_{\sharp} X, \ \ \text{where}\ \ f_{\sharp}(X)=f_{!}(X\sma S^A).$$
That is, for an $N$-free $G$-spectrum $X$ indexed on $f^*\sV^N$, $(i^*i_*X)^N\iso \SI^A X/N$. 
\end{thm}

It is usual to write this in the equivalent form $(i^*\SI^{-A}i_*X)^N\iso X/N$, but the
present form is more convenient for applications and more sensible from the categorical 
point of view. This looks enough like the formal Wirthm\"uller isomorphism to expect a 
similar proof. However, I do not have a helpful formal analysis.

\section{Change of universe}

Let $i\colon\sV\rtarr \sU$ be a map of $G$-universes, say for
definiteness an inclusion. We have the adjoint pair $(i_*,i^*)$
relating the stable homotopy categories $\sC$ and $\sD$ of 
$G$-spectra indexed on $\sV$ and $G$-spectra indexed on $\sU$. 
The left adjoint $i_*$ preserves compact objects, hence the
right adjoint $i^*$ preserves coproducts \cite[5.4]{May}.
Since $\sD$ is compactly generated, the right adjoint $i^*$
has a right adjoint $i_{!}$ by the triangulated adjoint functor
theorem \cite[6.3]{May}. It occurred to Po Hu to study
the structure of such functors $i_{!}$, and the natural question 
to ask is whether or not the formal Grothendieck isomorphism
theorem \cite[6.4]{May} applies. The functors 
$(i_*,i^*,i_{!})$ here play the roles of the functors 
$(f^*,f_*,f^{!})$ there, and the projection
formula takes the form of a natural isomorphism
$$Y\sma i^*X \iso i^*(i_*Y\sma X).$$
This holds when $Y$ is a suspension $G$-spectrum because the 
suspension $G$-spectrum functors on the two universes satisfy 
$i_*\SI^{\infty}\iso \SI^{\infty}$ and $i^*$ commutes with smash
products. It therefore holds in general \cite[5.6]{May}.

There is a natural map 
$$\ph: i_*Y\sma i_{!}Z \rtarr i_{!}(Y\sma Z)$$ 
which is an isomorphism for all dualizable $Y$ \cite[3.9]{May}. We ask whether or
not it is an isomorphism for all $Y$, and the answer is no. 
Indeed, the necessary hypothesis that $i_{!}$ preserves coproducts 
in \cite[6.4]{May} is satisfied if and only if $i^*$ takes 
generators to compact objects \cite[5.4]{May}, and this 
fails in general. To see this, let $\sU$ be a complete $G$-universe
and $\sV = \sU^G$ be the $G$-fixed subuniverse $\{U^G\}$, 
$U\in\sU$. The $G$-fixed point functor from $\sD$ to the stable
homotopy category of spectra is the composite of $i^*$ and the 
$G$-fixed point spectrum functor from $\sC$ to spectra. The latter 
functor preserves compact objects by inspection. For example,
it commutes with the suspension spectrum functor and therefore
takes suspension spectra of compact $G$-spaces to suspension
spectra of compact spaces. However, for a based $G$-space
$Y$, the $G$-spectrum $\SI^{\infty}Y$ indexed on $\sU$ has 
$G$-fixed point spectrum the wedge over conjugacy classes $(H)$
of the suspension spectra of the spaces 
$EWH_+\sma_{WH}\SI^{Ad(WH)}Y^H$ \cite[V.11.1]{LMS}. 
Even when $Y=S^0$, this spectrum is not compact. Therefore
$i^*\SI^{\infty}Y$ cannot be a compact $G$-spectrum indexed 
on $\sU^G$.

\section{Completion of the proof of the Wirthm\"uller isomorphism}

We must verify the hypotheses of \cite[4.14]{May} for $X= \SI^{\infty}H/K_+$. 
This means that, with $f_{\sharp}X = G_+\sma_H (X\sma S^{-L})$, we must construct a map
$\xi: f^*f_{\sharp} X \rtarr X$ such that certain diagrams commute. We need some space 
level constructions from \cite{LMS} to do this. The tubular neighborhood (\ref{tube}) 
gives a Pontryagin-Thom $G$-map $t: S^V\rtarr G_+\sma_H S^W$. It is $V$-dual to the 
counit $G$-map $\si: f_{!}f^*S^0 = G/H_+ \rtarr S^0$ \cite[III.5.2]{LMS}. The following 
construction, which is \cite[II.5.5]{LMS}, specializes to
give the $V$-dual $u: G_+\sma_H S^W\rtarr S^V$ to the unit 
$H$-map $\ze: S^0\rtarr f^*f_{!}S^0 = G/H_+$. We omit $f^*$ 
from notations, instead stating the equivariance explicitly.

\begin{con}
Let $H\times H$ act on $G$ by $(h_1,h_2)g = h_1gh_2^{-1}$ and 
act on $L\times H$ by $(h_1,h_2)(\la,h) = (h_1\la,h_1hh_2^{-1})$.
We think of the first factor $H$ as acting from the left, the 
second as acting from the right. Using the exponential map, 
construct an embedding $j: L\rtarr G$ of $L$ as a slice at $e$
such that 
\begin{equation}\label{jprop}
j(h\la) =hj(\la)h^{-1}\ \ \ \text{and}\ \ \ j(-\la) = j(\la)^{-1}.
\end{equation} 
Define $\tilde{j}: L\times H\rtarr G$ by 
$\tilde{j}(\la,h)=j(\la)h$. Then $\tilde{j}$ is an 
$(H\times H)$-map that embeds $L\times H$ onto an open neighborhood
of $e$. Collapsing the complement to a point, we obtain an 
$(H\times H)$-map $u: G_+\rtarr S^L\sma H_+$. For a based $H$-space $X$,
we obtain an induced (left) $H$-map
$$u: G_+\sma_H X \rtarr (S^L\sma H_+)\sma_H X \iso S^L\sma X.$$
Setting $X=S^W$ and identifying $S^V$ with $S^L\sma S^W$, we obtain
the promised $V$-dual $u:G_+\sma_H S^W\rtarr S^V$ of $\ze: S^0\rtarr G/H_+$.  
\end{con}

We only need the following definition and lemmas for the $H$-spaces
$X = H/K_+$, but it is simpler notationally to proceed more generally.
We write suspension coordinates on the right, $\SI^V Y = Y\sma S^V$. 
This is important to remember for control of signs, which are units in
Burnside rings. Recall that, for a
$G$-space $Y$ and an $H$-space $Z$, the natural isomorphism of $G$-spaces
$$\bar{\pi}: G_+\sma_H (Y\sma Z)\rtarr Y\sma (G_+\sma_H Z)$$
is given by the formulas
\begin{equation}\label{alpha} 
\bar{\pi}(g\sma y\sma z)=gy\sma g\sma z \ \ \ \text{and}\ \ \ 
\bar{\pi}^{-1}(y\sma g\sma z)=g\sma g^{-1}y\sma z.
\end{equation}

\begin{defn}\label{xiX} For a based $H$-CW complex $X$, define an $H$-map
$$\xi: f^*f_{\sharp} \SI^{\infty}_H X\rtarr \SI^{\infty}_H X$$
as follows. Observe that we have natural isomorphisms of $H$-spectra
$$
\SI^V f^*f_{\sharp}\SI^{\infty}_H X 
\iso S^V\sma (G_+\sma_H(X\sma S^{-L})) \hspace{1.1in}$$
$$\hspace*{1.1in} \overto{\bar{\pi}^{-1}} 
 G_+\sma_H(S^V\sma X\sma S^{-L}) \iso \SI^{\infty}_H (G_+\sma_H(X\sma S^W))$$
and 
$$
\SI^V\SI^{\infty}_H X\iso \SI^{\infty}_H(X\sma S^V).$$
Let $\SI^V\xi$ correspond under these isomorphisms to $\SI^{\infty}_Hu$, where
$$u: G_+\sma_H (X\sma S^W)\rtarr S^L\sma X\sma S^W\iso X\sma S^V.$$
\end{defn}

The following observation is taken from \cite[II.5.8]{LMS}. 

\begin{lem} For a $G$-space $Y$ and an $H$-space $Z$, the $H$-map 
$$\xymatrix@1{
G_+\sma_H (Y\sma Z) \iso Y\sma (G_+\sma_H Z)\ar[r]^-{\id\sma u} & Y\sma S^L\sma Z}$$
is canonically $H$-homotopic to the $H$-map 
$$\xymatrix@1{G_+\sma_H (Y\sma Z)\ar[r]^-u  & S^L\sma Y\sma Z\iso Y\sma S^L\sma Z.}$$
Therefore, for any $H$-map $\tha:Y\rtarr X$ from a $G$-space $Y$ to an $H$-space $X$,
the following diagram is canonically homotopy commutative. 
$$\xymatrix{
G_+\sma_H (Y\sma Z)\ar[d]_{\id\sma(\tha\sma\id)} \ar@{} [r] |{\iso} 
& Y\sma (G_+\sma_H Z)\ar[r]^-{\id\sma u} & Y\sma S^L\sma Z \ar[d]^{\tha\sma\id}\\
G_+\sma_H (X\sma Z) \ar[r]_{u} & S^L\sma X\sma Z \ar @{} [r] |{\iso} & X\sma S^L\sma Z}$$
\end{lem}
\begin{proof}
Both maps send all points of $G$ not in $\tilde{j}(L\times H)$ to the basepoint.
The first takes $\tilde{j}(\la,h)\sma y\sma z$ to $j(\la)hy\sma \la\sma hz$
and the second takes it to $hy\sma \la\sma hz$.
Applying $j$ to a contracting homotopy of $L$, we obtain an $H$-homotopy
from $j:L\rtarr G$ to the constant map at $e$, giving the required homotopy.
The last statement follows since the diagram commutes by naturality if its top 
arrow is replaced by $u$.
\end{proof}

This leads to the following naturality statement. The partial naturality 
diagram \cite[4.16]{May} is the case in which $Y=f_*X$ and $\tha$ 
is the counit of the $(f^*,f_*)$ adjunction.

\begin{lem} Let $Y$ be a $G$-spectrum and $X$ be a based $H$-space.
For any map $\tha: Y\rtarr \SI^{\infty}_H X$ of $H$-spectra, the 
following diagram of $H$-spectra commutes in $\sC$.
$$
\xymatrix{
G_+\sma_H (Y\sma S^{-L}) \ar[r]^-{\xi} \ar[d]_{\id\sma(\tha\sma\id)} & Y \ar[d]^{\tha} \\
G_+\sma_H (\SI^{\infty}_HX\sma S^{-L}) \ar[r]_-{\xi} 
& \SI^{\infty}_HX \\}
$$
\end{lem}
\begin{proof} The upper map $\xi$ in the diagram is defined formally in \cite[4.9]{May}.
It suffices to prove that the diagram commutes after suspension by $V$. This has the effect 
of replacing
$S^{-L}$ by $S^W$ on the left and suspending by $V$ on the right. Comparing 
\cite[4.9]{May} with Definition \ref{xiX} and taking $Z=S^W$, the conclusion reduces to 
application of the spacewise diagram of the previous lemma to the spaces that comprise
the given spectra. Technically, this is most easily seen using prespectra or orthogonal
spectra rather than actual spectra, but the essential point is just that the homotopy
in the previous lemma is sufficiently natural to commute with the structure maps.
\end{proof} 

To complete the proof of Theorem \ref{Wirth}, it suffices to show that the following 
specialization of the diagram \cite[4.17]{May} commutes.
\begin{equation}\label{dia}
\xymatrix{
\SI^{\infty}_H X\sma S^{-L} \ar[r]^-{\ze} \ar[d]_{\ze} 
& f^*f_{!}(\SI^{\infty}_H X\sma S^{-L}) \\
f^*f_{!}(\SI^{\infty}_H X\sma S^{-L}) \ar[r]_-{f^*\ta} 
& f^*f_{!}(f^*f_{!}(\SI^{\infty}_H X\sma S^{-L})\sma S^{-L}) 
\ar[u]_{f^*f_{!}(\xi\sma\id)} \\}
\end{equation}
Here $f_{!}(-) = G_+\sma_H(-)$, $f^*$ is the forgetful functor, $\ze$ is the unit of 
the $(f_{!},f^*)$ adjunction, and $\ta$ is the map defined in \cite[4.6]{May} with 
$Y=f_{!}(\SI^{\infty}_H X\sma S^{-L})$. We shall see that this reduces 
to the following space level observations from \cite[II\S5]{LMS}.

\begin{lem}\label{ret} The following composite is $H$-homotopic to the identity map.
$$\xymatrix@1{ 
S^V \ar[r]^-{t} & G_+\sma_H S^W \ar[r]^-{u} & S^V.\\}$$
\end{lem} 
\begin{proof}
Composing the embeddings $\tilde{i}: G\times_H W\rtarr V$ and 
$\tilde{j}: L\times H\rtarr G$, we obtain an embedding
$V = L\times W =(L\times H)\times _H W \rtarr G\times_H  V$.
The composite $u\com t$ is $k^{-1}$ on $k(V)$, and it collapses the complement
of $k(V)$ to the basepoint. The embedding $k$ is isotopic to the identity, 
and application of the Pontryagin-Thom construction to an isotopy gives a
homotopy $\id\htp u\com t$. 
\end{proof}

\begin{lem}\label{muck} 
For an $H$-space $X$, the following diagram is $H$-homotopy commutative.
Here $\si:S^L\rtarr S^L$ maps $\la$ to $-\la$.
\small{
$$\xymatrix{
\SI^VX \ar[r]^-{\SI^V\ze} \ar[dd]_{\SI^V\ze} & \SI^V(G_+\sma_H X) \iso S^V\sma(G_+\sma_H X)
& G_+\sma_H (S^V\sma X) \ar[l]_-{\bar{\pi}^{-1}} \\
& & G_+\sma _H \SI^W(S^L\sma X) \ar[u]_{\iso}\\
\SI^V(G_+\sma_H X) \ar[r]_-{\id\sma t}  & 
(G_+\sma_H X)\sma (G_+\sma_H S^W) \ar[r]_-{\bar{\pi}^{-1}} 
& G_+\sma_H \SI^W(G_+\sma_H X) \ar[u]_{\id\sma\SI^{W}(\si\sma\id)u}\\}$$
}
\end{lem}
\begin{proof}
The composite around the bottom maps all points with $V$ coordinate not in $\tilde{i}(j(L)\times W))$
to the basepoint. It maps the point $x\sma \tilde{i}(j(\la), w)$ to $(j(\la,x)\sma j(\la)(\la,w)$; 
the sign map $\si$ enters due to (\ref{jprop}) and (\ref{alpha}). The $H$-contractibility of $L$
implies that an $H$-homotopic map is obtained if we replace $f(\la)$ by the identity element $e\in G$.
Thus the composite around the right is $H$-homotopic to $\ze\sma u\com t$, which is $H$-homotopic
to the identify by Lemma \ref{ret}.
\end{proof}

\begin{proof}[Proof of Theorem \ref{Wirth}] We must show that the diagram (\ref{dia}) commutes.
It suffices to prove this after suspending by $V$ and 
replacing $X$ by $\SI^VX$. This has the effect of replacing the $H$-spectra $S^{-L}$ that 
appear in the diagram by the $H$-spaces $S^W$. Since the functors appearing in the diagram 
commute with the respective suspension spectrum functors, this reduces the problem to the 
space level. Here a slightly finicky diagram chase, which amounts to a check of signs
coming from permutations of suspension coordinates, shows that the resulting diagram
commutes by Lemma \ref{muck}. One point is
that $\si\sma\id: S^L\sma S^L\rtarr S^L\sma S^L$ is $H$-homotopic to the transposition 
via the homotopy given by multiplying by the $(2\times 2)$-matrices with rows $(-t,1-t)$
and $(1-t,t)$ for $t\in I$. Another is that we must apply Lemma \ref{muck} with $X$ replaced
by $X\sma S^W$, which has the effect of introducing a permutation of $S^L$ past $S^W$. 
In more detail, after applying $\SI^{\infty}_H$ and replacing $X$ by $X\sma S^W$, the bottom 
composite in the diagram of the  previous lemma agrees with 
$$\SI^V\ta: \SI^V(G_+\sma_H (\SI^V X\sma S^{-L})) 
\rtarr \SI^V(G_+\sma_H(G_+\sma_H(\SI^VX\sma S^{-L})\sma S^{-L})$$
under the evident isomorphisms of its domain and target. The interpretation of the right
vertical composite is trickier, because of the specification of 
$$\xi: G_+\sma_H (\SI^V X\sma S^{-L})\rtarr \SI^V X$$ 
in Definition \ref{xiX}. A diagram  chase after suspending by $V$ shows that $\xi$
agrees under the evident isomorphism of its source with
$$(\si\sma\id)\com u: G_+\sma_H (X\sma S^W)\rtarr S^L\sma X\sma S^W \iso \SI^V X.$$ 
Notice that we have evident isomorphisms
$$\SI^{\infty}_H\SI^V(G_+\sma_H (X\sma S^W)) \iso \SI^V(G_+\sma_H(\SI^VX\sma S^{-L}))
\iso \SI^{\infty}_HG_+\sma_H (\SI^V X\sma S^W),$$
the first of which internally expands $S^W$ to $S^V\sma S^{-L}$ and the second of which 
leaves $\SI^VX$ alone but uses $\bar{\pi}$ to bring $S^V$ inside and contracts 
$S^{-L}\sma S^V$ to $S^W$.
Their composite enters into the upper right corner of the required diagram chase; the transpositions 
that appear cancel out others, resulting in a sign free conclusion. With these indications, the 
rest is routine.
\end{proof}

\begin{rem} Effectively, this proof uses the space level arguments of \cite[II\S5]{LMS}, but
eliminates the need for the spectrum level arguments of \cite[II\S6]{LMS}. We warn the reader 
that \cite[II.5.2]{LMS} and hence the first diagram of \cite[II.6.12]{LMS} are incorrect.
Fortunately, they are also irrelevant. 
\end{rem}

\end{document}